\documentclass{article}
\usepackage{xcolor}
\usepackage{hyperref}
\usepackage{amsmath, amsfonts, amssymb,amsthm}
\usepackage{mathrsfs}
\usepackage[margin = 1 in]{geometry}
\usepackage{float}
\restylefloat{table}
\usepackage{changepage}
\usepackage{lipsum}

\newcommand{\Hom}{\mathrm{Hom}}
\newcommand{\Proof}{\noindent \textbf{Proof:}}

\renewcommand{\span}{\textrm{span}}

\newcommand{\fdlcat}[1]{\textrm{Rel}(#1)}
\newcommand{\REL}{\textrm{Rel}}

\theoremstyle{definition}
\newtheorem{defn}{Definition}[section]

\theoremstyle{plain}
\newtheorem{lemma}{Lemma}[section]

\newtheorem{thm}{Theorem}[section]
\newtheorem{prop}{Proposition}[section]

\newtheorem{cor}{Corollary}[section]

\makeatletter
\newenvironment{amatrix}{%
  \left\langle
  \env@matrix
}{%
  \endmatrix
  \right\rangle
}

\newcommand{\amatrixarrow}[3][.5ex]{%
  \def\amatrixarrow@shift{#1}%
  \def\amatrixarrow@left{#2}%
  \def\amatrixarrow@right{#3}%
  \collect@body\amatrixarrow@next
}
\newcommand*{\amatrixarrow@next}[1]{%
  \mathpalette{\amatrixarrow@}{#1}%
}
\newcommand*{\amatrixarrow@}[2]{%
  \sbox0{$\m@th#1\begin{amatrix}#2\end{amatrix}$}%
  \sbox2{$\m@th#1\begin{matrix}#2\end{matrix}$}%
  \sbox4{$\m@th
    \amatrixarrow@style{#1}%
    \amatrixarrow@left
  $}%
  \sbox6{$\m@th
    \amatrixarrow@style{#1}%
    \amatrixarrow@right
  $}%
  \sbox8{$\m@th\amatrixarrow@style{#1}\text{\kern\amatrixarrow@shift}$}%
  \dimen0=.5\dimexpr\wd0-\wd2\relax
  \vtop{%
    \hbox{%
      \kern.5\dimexpr\wd4-\dimen0\relax
      \copy0 %
    }%
    \kern\wd8 %
    \nointerlineskip
    \hbox to \dimexpr\wd0+.5\wd4+.5\wd6-\dimen0\relax{%
      \copy4 %
      \hfill
      \sbox2{$\m@th
        \amatrixarrow@style{#1}%
        {}\xrightarrow{}{}%
      $}
      $\m@th
        \amatrixarrow@style{#1}%
        \xrightarrow{%
          \kern\dimexpr\wd0-.5\wd4-.5\wd6-\dimen0-\wd2\relax
        }%
      $%
      \hfill
      \copy6 %
    }%
  }%
}
\newcommand*{\amatrixarrow@style}[1]{%
  \ifx#1\displaystyle
    \scriptstyle
  \else\ifx#1\textstyle
    \scriptstyle
  \else
    \scriptscriptstyle
  \fi\fi
}
\makeatother

\begin{document}

\begin{center}
{\Large {\bfseries New Examples of Dimension Zero Categories}} \\
Andrew Gitlin \\
March 2018 \\
\end{center}

\begin{center}
\textbf{Abstract}
{\small 
\begin{adjustwidth}{1in}{1in}
We say that a category $\mathscr{D}$ is dimension zero over a field $F$ provided that every finitely generated representation of $\mathscr{D}$ over $F$ is finite length.  We show that $\fdlcat{R}$, a category that arises naturally from a finite idempotent semiring $R$, is dimension zero over any infinite field.  One special case of this result is that $\REL$, the category of finite sets with relations, is dimension zero over any infinite field.  
\end{adjustwidth}}
\end{center}

\section{Introduction and Preliminaries}

\indent \indent We define a representation of a category $\mathscr{D}$ over a field $F$ to be a functor from $\mathscr{D}$ to $Vect_{F}$, the category of vector spaces over $F$.  We say that a category $\mathscr{D}$ is dimension zero over a field $F$ provided that every finitely generated representation of $\mathscr{D}$ over $F$ is finite length.  The purpose of this paper is to show that $\REL$, the category of finite sets with relations, is dimension zero over any infinite field.  Our method of argument allows this result to be generalized to categories that we call $\fdlcat{R}$, where $R$ is any finite idempotent semiring (Definition 1.8).  Bouc and Th\'{e}venaz \cite{correspond} have independently shown that $\REL$, the category of finite sets with relations, is dimension zero over any field.  Theorem 3.2 states that $\fdlcat{R}$ is dimension zero over any infinite field for any finite idempotent semiring $R$. 

\par
For the rest of this paper, let $\mathscr{D}$ be a combinatorial category, i.e. a category such that $\Hom(a,b)$ is finite for all objects $a,b \in \mathscr{D}$, and let $F$ be a field.  We will let $Vect_F$ denote the category of vector spaces over $F$; the objects are vector spaces over $F$ and the morphisms are linear transformations.  Finally, for the rest of this paper, let $[n]$ be the set $\{1,...,n\}$ for any whole number $n$.

\par
We will now introduce several notions in representation theory which will be important in this paper.

\begin{defn}
	A \underline{representation} of $\mathscr{D}$ over $F$ is a functor from $\mathscr{D}$ to $Vect_{F}$. 
\end{defn}

\par Concretely, a representation $V$ of $\mathscr{D}$ over $F$ takes every object $d \in \mathscr{D}$ to a vector space $V(d)$ over $F$, takes every morphism $g \in \Hom(d,e)$ to a linear map $V(g) \in \Hom(V(d),V(e))$ for all $d,e \in \mathscr{D}$, and satisfies the following two properties.
\begin{itemize}
	\item $V(f \circ g) = V(f) \circ V(g)$ \textrm{for all $f \in \Hom(y,z), g \in \Hom(x,y)$ for all $x,y,z \in \mathscr{D}$}
	\item $V(Id_d) = Id_{V(d)}$ \textrm{for all $d \in \mathscr{D}$}
\end{itemize}

\begin{defn}
	Let $V$ be a representation of $\mathscr{D}$ over $F$.  A \underline{subrepresentation} of $V$ is a representation $W$ of $\mathscr{D}$ over $F$ such that $W(d)$ is a vector subspace of $V(d)$ for all $d\in\mathscr{D}$ and $W(f)$ is the restriction of $V(f)$ to $W(d)$ for all $d,d'\in\mathscr{D}$ and $f\in\Hom_{\mathscr{D}}(d,d')$.
\end{defn}

\par
Two particularly easy examples of a representation of $\mathscr{D}$ over $F$ are the zero representation and the trivial representation.  The zero representation sends every object of $\mathscr{D}$ to $0$ and every morphism in $\mathscr{D}$ to the zero transformation.  The trivial representation sends every object of $\mathscr{D}$ to $F$ and every morphism in $\mathscr{D}$ to the identity transformation.

\begin{defn}
	A representation $V$ is \underline{irreducible} provided that $V$ is not the zero representation and that the only subrepresentations of $V$ are the zero representation and $V$ itself.
\end{defn}

\begin{defn}
	A representation $V$ of $\mathscr{D}$ over $F$ is \underline{finitely generated} provided that there exist objects $d_1,...,d_i \in\mathscr{D}$ and $v_{1,1},...,v_{1,j_1} \in V(d_1),...,v_{i,1},...,v_{i,j_i} \in V(d_{i})$ such that if $W$ is a subrepresentation of $V$ and $v_{1,1},...,v_{1,j_1} \in W(d_1),...,v_{i,1},...,v_{i,j_i} \in W(d_{i})$ then $W = V$.
\end{defn}

\begin{defn} 
	A representation $V$ is \underline{finite length} provided that any non-repetitive chain of subrepresentations of $V$ is finite.
\end{defn}

\par
An equivalent definition of finite length is that a representation $V$ is finite length provided that there exists some non-repetitive finite chain $0=W_{0}\subsetneq...\subsetneq W_{n}=V$ of subrepresentations of $V$ such that each $W_{i+1} / W_{i}$ is irreducible.  When this is the case, $W_0,...,W_n$ is called a composition series for $V$ and $n$ is called the length of $V$.  The Jordan-H\"older Theorem guarantees that if $W'_0,...,W'_m$ is another composition series for $V$ then $m = n$ (and thus the length of $V$ is well-defined) and in fact that the $W'_{i+1} / W'_{i}$ are a permutation of the $W_{i+1} / W_{i}$.  A statement and proof of the Jordan-H\"older Theorem, in a more general setting, can be found as Theorem 2.1 in $\cite{JH}$; since $Vect_F$ is an abelian category, the category of functors from $\mathscr{D}$ to $Vect_F$ is also an abelian category and thus this theorem applies here.

\par
\begin{defn}
	A category $\mathscr{D}$ is \underline{dimension zero} over a field $F$ provided that every finitely generated representation of $\mathscr{D}$ over $F$ is finite length.
\end{defn}

\par
This paper considers a natural construction of a category from a finite idempotent semiring.  We will now define what a finite idempotent semiring is and explain this construction.

\begin{defn}
	A \underline{finite idempotent semiring} is a finite set $R$ equipped with two binary operations, denoted $+$ (addition) and $*$ (multiplication), satisfying the following axioms.
	\begin{itemize}
		\setlength\itemsep{0em}
		\item $+$ is commutative, associative, and idempotent and there exists an additive identity $0 \in R$
		\item $*$ is associative and there exists a multiplicative identity $1 \in R$
		\item $*$ distributes over $+$
		\item $0*a = a*0 = 0$ for all $a \in R$
	\end{itemize}
\end{defn}

\par
For the rest of this paper, let $R$ be a finite idempotent semiring.  

\begin{defn}
We will now define a category, which we will denote $\underline{\fdlcat{R}}$, which arises naturally from the finite idempotent semiring $R$.  The objects are the whole numbers.  For any whole numbers $x,y$, a morphism from $x$ to $y$ is a $x \times y$ matrix with elements of $R$ as its entries.  The composition of morphisms is given by matrix multiplication.  
\end{defn}

\par
Throughout the rest of this paper, for any whole numbers $x,y$, if $A \in \Hom_{\fdlcat{R}}(x,y)$ and the $(i,j)$ entry of $A$ is $a_{i,j}$ for all $i \in [x],j \in [y]$, then we will let $(a_{i,j})$ denote the morphism $A$.  Using this notation, we can now express the rule for composing morphisms in $\fdlcat{R}$ more explicitly.  If $x,y,z$ are whole numbers and $A = (a_{i,j}) \in \Hom(x,y)$ and $B = (b_{i,j}) \in \Hom(y,z)$, then
\begin{equation*}
B \circ A = AB = \left ( \sum_{k = 1}^y a_{i,k} * b_{k,j} \right ) \in \Hom(x,z).
\end{equation*}

\par
A special case of the above discussion is when $R = \{0,1\}$ and $+$ and $*$ are given by logical OR and logical AND, respectively.  In this special case, $\fdlcat{R}$ is called the category of finite sets with relations and is denoted $\REL$.

\par 
Another example of a finite idempotent semiring is the truncated tropical semiring $R = \{0,1,...,n,\infty\}$ where $n$ is a fixed whole number.  Addition $\oplus$ and multiplication $\otimes$ on $R$ are defined as follows.
\begin{equation*}
x \oplus y := \min(x,y) \indent x \otimes y := 
\left\{
\begin{array}{ll}
      \min(x+y,n) & \textrm{if } x,y \neq \infty \\
      \infty & \textrm{if } x = \infty \textrm{ or } y = \infty
\end{array}
 \right. 
\end{equation*}
\noindent The truncated tropical semiring is a truncated version of the tropical semiring $\mathbb{R} \cup \{\infty\}$ where the addition operation $\oplus$ is given by $x \oplus y := \min(x,y)$ and the multiplication operation $\otimes$ is given by $x \otimes y := x+y$.  The reason for the truncation is that a finite idempotent semiring must be a finite set.  For an introduction to the tropical semiring, see $\cite{speyer}$. 

\pagebreak

\section{A Partial Order on $\fdlcat{R}$}

\begin{defn}
For any $a,b \in R$, we write $a \subseteq b$ and $b \supseteq a$ when $a+b = b$.
\end{defn}

\indent The following lemma lists some important properties of $\subseteq$. 

\begin{lemma}
\hfill 
	\begin{enumerate}
	\item $\subseteq$ is a partial order on $R$ with minimal element $0$
	\item for all $a,b \in R$, we have $a \subseteq a + b$
	\item for all $a,b,c \in R$, if $a,b \subseteq c$ then $a+b \subseteq c$
	\end{enumerate}
\end{lemma}

\indent The proof of Lemma 2.1 is a series of routine computations and is left as an exercise for the reader.

\section{$\fdlcat{R}$ is Dimension Zero over any Infinite Field}

\par
\indent \indent Fix $d,x,y\in\mathscr{D}$.  For any $s\in\Hom(x,x)$,  define a matrix $M_{s}$ with rows and columns indexed by $\Hom(d,x)$ by letting the $(f,g)$ entry be $1$ if $s\circ f=g$ and $0$ otherwise for all $f,g\in\Hom(d,x)$.  We write $x\leq_{d} y$ when $\span_{F}\{M_{t} : t\in\Hom(x,y,x)\}$ contains the identity matrix, where $\Hom(x,y,x)$ is defined to be $\{\lambda\in\Hom(x,x) :$ there exists $a\in\Hom(x,y), b\in\Hom(y,x)$ such that $\lambda = b\circ a\}$.  

\par
Proposition 2.5 in \cite{WG15} says that $\leq_{d}$ is a preorder on the objects of $\mathscr{D}$.  Furthermore, Theorem 1.2 in \cite{WG15} uses $\leq_{d}$ to provide a criterion for determining whether or not $\mathscr{D}$ is dimension zero over $F$.

\begin{thm} {\textup{(\cite{WG15}, Thm. 1.2)}}
	A category $\mathscr{D}$ is dimension zero over a field $F$ if and only if $\Hom_{\mathscr{D}}(a,b)$ is finite for all $a,b\in\mathscr{D}$ and for all $d\in\mathscr{D}$ there exists a finite set $Y_{d}$ of objects of $\mathscr{D}$ such that for all $x\in\mathscr{D}$ there exists some $y\in Y_{d}$ such that $x\leq_{d} y$ over $F$.
\end{thm}

\par
Theorem 1.2 in \cite{WG15} was proven in the more general setting of representations of categories over rings.

\par 
Proposition 2.1 in \cite{WG15} gives one useful property of $\leq_{d}$.

\begin{prop} {\textup{(\cite{WG15}, Prop. 2.1)}}
Let $d,x,y\in\mathscr{D}$.  If $\Hom(d, x)$ is finite, then $\span_{F}\{M_{t}$ : $t\in\Hom(x,y,x)\}$ contains an invertible matrix if and only if $x \leq_{d} y$.
\end{prop}

\par
The following proposition provides one method for proving that $x \leq_{d} y$ for some fixed $d,x,y \in \mathscr{D}$.

\begin{prop}
	Let $d,x,y\in\mathscr{D}$.  If $\Hom(d,x)$ is finite and there exists a partial order $\preceq$ on $\Hom(d,x)$ and a function $s: \Hom(d,x) \rightarrow \Hom(x,y,x)$ such that $s(f) \circ f = f$ and $s(f) \circ h \succeq h$ for all $f,h \in\Hom(d,x)$, then $x\leq_{d}y$ over any infinite field $A$.
\end{prop}

\Proof 

{

Extend $\preceq$ to a total order $\leqslant$ on $\Hom(d,x)$.  For all $F \in \Hom(x,x)$, let the rows and columns of $M_F$ be arranged from least to greatest according to $\leqslant$.  For all $f,g \in \Hom(d,x)$, let $b_{f,g}$ be the $(g,g)$ entry of $M_{s(f)}$.  Note that $b_{f,f} = 1$ for all $f \in \Hom(d,x)$ since $s(f) \circ f = f$ for all $f \in \Hom(d,x)$.  Also note that $b_{f,g}$ is either 1 or 0 for all $f,g \in \Hom(d,x)$.

We will first show that there exists $\{a_f \in A : f \in \Hom(d,x)\}$ such that 
\begin{equation*}
\sum \limits_{f \in \Hom(d,x)} b_{f,g}a_f \neq 0
\end{equation*}
for all $g \in \Hom(d,x)$.  Let $m = |\Hom(d,x)|$ and let $\Hom(d,x) = \{f_1,...,f_m\}$.  It is enough to show that for all $n \in [m]$ there exist $a_{f_1},...,a_{f_n} \in A$ such that 
\begin{equation*}
\sum \limits_{i = 1}^{n} b_{f_i,f_1}a_{f_i} \neq 0,...,\sum \limits_{i = 1}^{n} b_{f_i,f_{n}}a_{f_i} \neq 0.
\end{equation*}
\noindent We will use induction on $n$.  If $n = 1$, then, since $b_{f_1,f_1} = 1$, we have that $a_{f_1} = 1$ is a solution to $b_{f_1,f_1}a_{f_1} \neq 0$ as desired.  Suppose that $n \geq 2$ and the result holds for $n-1$.  By the inductive hypothesis, there is a solution $a_{f_1},...,a_{f_{n-1}}$ to the system of equations
\begin{equation*}
\sum \limits_{i = 1}^{n-1} b_{f_i,f_1}a_{f_i} \neq 0,...,\sum \limits_{i = 1}^{n-1} b_{f_i,f_{n-1}}a_{f_i} \neq 0.
\end{equation*}
\noindent For all $k \in [n]$, let $S_k = - \sum_{i = 1}^{n-1} b_{f_i,f_k}a_{f_i}$.  Note that $S_k \neq 0$ for all $k \in [n-1]$.  Since $A$ is infinite and $\{S_k : k \in [n]\}$ is finite, there exists $a_{f_n} \in A \backslash \{S_k : k \in [n]\}$.  For all $k \in [n-1]$, we have that
\begin{equation*}
\begin{array}{ll}
&\sum \limits_{i = 1}^{n} b_{f_i,f_k}a_{f_i} = b_{f_n,f_k}a_{f_n} + \sum \limits_{i = 1}^{n-1} b_{f_i,f_k}a_{f_i} = b_{f_n,f_k}a_{f_n} - S_k =
\left \{
\begin{array}{ll}
	-S_k & \textrm{if } b_{f_n,f_k} = 0 \\
  	a_{f_n} - S_k & \textrm{if } b_{f_n,f_k} = 1
\end{array}
\neq 0.
\right.
\end{array}
\end{equation*}
\noindent Furthermore, noting that $b_{f_n,f_n} = 1$, we have that
\begin{equation*}
\sum \limits_{i = 1}^{n} b_{f_i,f_n}a_{f_i} = b_{f_n,f_n}a_{f_n} + \sum \limits_{i = 1}^{n-1} b_{f_i,f_n}a_{f_i} = a_{f_n} - S_n \neq 0.
\end{equation*}
\noindent Thus, as desired, $a_{f_1},...,a_{f_n}$ is a solution to the system of equations $\sum_{i = 1}^{n} b_{f_i,f_1}a_{f_i} \neq 0,...,\sum_{i = 1}^{n} b_{f_i,f_{n}}a_{f_i} \neq 0$.

Let $X = \sum_{f\in\Hom(d,x)} a_f M_{s(f)}$.  For all $g \in \Hom(d,x)$, the $(g,g)$ entry of $X$ is $\sum_{f \in \Hom(d,x)} b_{f,g}a_f$, which is non-zero.  Furthermore, since $s(f) \circ h \succeq h$ for all $f,h \in\Hom(d,x)$, $M_{s(f)}$ is upper triangular for all $f \in \Hom(d,x)$ and thus $X$ is upper triangular.  Thus, since $X$ is an upper triangular matrix with all of its diagonal entries being non-zero, $X$ is invertible.  Therefore, since $X\in \span_A\{M_{t}$ : $t\in\Hom(x,y,x)\}$, we are done by Proposition 3.1. $\square$ \\

}

\par
We are now ready to prove the following theorem, which is the main result of this paper.

\begin{thm}
	If $R$ is a finite idempotent semiring, then $\fdlcat{R}$ is dimension zero over any infinite field.
\end{thm} 

\Proof {

\par
Let $R$ be a finite idempotent semiring.  

\par
For all $x,y \in \fdlcat{R}$, let $0_{x \times y}$ be the $x \times y$ matrix with each entry being 0.  For all $x \in \fdlcat{R}$, let $Id_x$ be the $x \times x$ matrix with each diagonal entry being 1 and every other entry being 0.  For all $x,y,z \in \fdlcat{R}$, $0_{x \times y}A = 0_{x \times z}$ for all $A \in \Hom(x,y)$ and $B0_{y \times z} = 0_{x \times z}$ for all $B \in \Hom(x,y)$.  For all $x,y \in \fdlcat{R}$, $Id_xA = A$ for all $A \in \Hom(x,y)$ and $BId_y = B$ for all $B \in \Hom(x,y)$.

\par
Let $n = |R|$.  It is enough to show $x\leq_{d}n^{d}$ for all $d,x\in \fdlcat{R}$, since then we can apply Theorem 3.1 with $Y_d = \{n^d\}$ for all $d \in \fdlcat{R}$.  Fix $d,x \in \fdlcat{R}$.  If $x\leq n^{d}$ then the below calculation shows that $Id_{x}\in\Hom(x,n^{d},x)$ and therefore, since $M_{Id_x}$ is the identity matrix, we are done. 
\begin{equation*}
\begin{bmatrix}
Id_x & 0_{x \times (n^d-x)} 
\end{bmatrix}
\begin{bmatrix}
Id_x \\ 0_{(n^d-x) \times x} 
\end{bmatrix}
= 
Id_x
\end{equation*}

\par
For the $x > n^{d}$ case, we will use Proposition 3.2.

\par
We will first construct $\preceq$.  For any $f = (a_{i,j}),g = (b_{i,j}) \in\Hom(d,x)$, we will write $f \preceq g$ if and only if $a_{p,q} \subseteq b_{p,q}$ for all $p \in [d],q \in [x]$.  Since $\subseteq$ is a partial order by Lemma 2.1, $\preceq$ is also a partial order.

\par
We will now construct $s : \Hom(d,x) \rightarrow \Hom(x,n^{d},x)$.  For any $f = (a_{i,j}) \in\Hom(d,x)$, let $s(f)$ be the matrix $M=(m_{i,j}) \in \Hom(x,x)$ where 
\begin{equation*}
m_{i,j}=
\left \{
\begin{tabular}{l}
	$1$ if $a_{k,i} \subseteq a_{k,j}$ for all $k\in[d]$ \\
  	$0$ otherwise
\end{tabular}
\right.
\end{equation*}
for $i,j\in[x]$.  The proposition below verifies that $s(f)\in\Hom(x,n^{d},x)$.

\pagebreak

\begin{prop}
	We have $s(f)\in\Hom(x,n^{d},x)$.
\end{prop}

\Proof {

\par
Since there are only $n^{d}$ distinct $d$-tuples with entries from $R$, there are at most $n^{d}$ distinct columns in $f$.  Note that if columns $p \in [x]$ and $q \in [x]$ of $f = (a_{i,j})$ are identical, then for all $i \in [x]$

\begin{equation*}
m_{i,p} = 
\left \{
\begin{tabular}{l}
	$1$ if $a_{k,i} \subseteq a_{k,p} = a_{k,q}$ for all $k\in[d]$ \\
  	$0$ otherwise
\end{tabular}
\right.
= m_{i,q} 
\end{equation*}

\noindent and thus columns $p$ and $q$ of $M$ are identical.  Therefore, $M$ has at most $n^{d}$ distinct columns.  Let $v$ be the number of distinct columns of $M$, noting that $v \leq n^d$.  Let $D \in \Hom(x,n^d)$ be the $x \times n^{d}$ matrix which has the $v$ distinct columns of $M$ as its first $v$ columns and has every entry in its $n^d - v$ remaining columns as $0$.  Let $E = (e_{i,j}) \in \Hom(n^d,x)$ be a $n^{d} \times x$ matrix where
\begin{equation*}
e_{i,j} = 
\left \{
\begin{tabular}{l}
	$1$ if $i \leq v$ and column $j$ of $M$ is column $i$ of $D$ \\
  	$0$ otherwise
\end{tabular}
\right.
\end{equation*}
\noindent for all $i \in [n^d],j \in [x]$.  For all $i,j \in [x]$, we have  
\begin{equation*}
\sum_{k \in [n^d]} d_{i,k} * e_{k,j} = \sum_{\substack{k \in [v] \\ \textrm{col. $j$ of $M$ is col. $k$ of $D$}}} d_{i,k} = \sum_{\substack{k \in [v] \\ \textrm{col. $j$ of $M$ is col. $k$ of $D$}}} m_{i,j} \stackrel{\dagger}{=} m_{i,j}.
\end{equation*}
\noindent The $\dagger$ step uses the idempotence of $+$.  Thus $M=DE$.  Thus $s(f)\in\Hom(x,n^{d},x)$ as desired. $\square$ \\

}

\par
We are left to show that, for all $A \in\Hom(d,x)$, we have $s(A)\circ A = A$ and $s(A)\circ B \succeq B$ for all $B\in\Hom(d,x)$.  Fix $A=(a_{i,j}) \in\Hom(d,x)$.  Let $M=(m_{i,j})$ be $s(A)$.  

\par
We will first show that $s(A) \circ B \succeq B$ for all $B\in\Hom(d,x)$.  Fix $B = (b_{i,j}) \in\Hom(d,x)$.  Let $C=(c_{i,j})$ be $s(A) \circ B$, noting that $C=BM$.  By definition, it is enough to show $c_{i,j} \supseteq b_{i,j}$ for all $i \in [d],j \in [x]$.  For all $j\in[x]$, we have $a_{k,j} \subseteq a_{k,j}$ for all $k\in[d]$ and thus $m_{j,j}=1$.  Thus, for all $i \in [d],j \in [x]$, we have $c_{i,j} = \sum_{l=1}^{x}(b_{i,l} * m_{l,j}) \supseteq b_{i,j} * m_{j,j}= b_{i,j} * 1=b_{i,j}$ as desired.

\par
We will now show that $s(A) \circ A = A$.  Let $C=(c_{i,j})$ be $s(A)\circ A$, noting that $C=AM$.  We have $s(A) \circ A \succeq A$ by the previous paragraph.  Therefore, it is enough to show that $s(A) \circ A \preceq A$.  By definition, it is enough to show that $c_{i,j} \subseteq a_{i,j}$ for all $i \in [d],j \in [x]$.  Fix $i \in [d],j \in [x]$.  Fix $k\in[x]$.  Note that $m_{k,j}$ is either $0$ or $1$.  If $m_{k,j} = 0$ then $a_{i,k} * m_{k,j} = 0 \subseteq a_{i,j}$.  If $m_{k,j} = 1$ then we have $a_{r,k} \subseteq a_{r,j}$ for all $r\in[d]$ and thus $a_{i,k} \subseteq a_{i,j}$ and thus $a_{i,k} * m_{k,j} = a_{i,k} \subseteq a_{i,j}$.  Therefore $a_{i,k} * m_{k,j} = a_{i,k} \subseteq a_{i,j}$ for all $k\in[x]$, independently of $m_{k,j}$.  Thus $c_{i,j} = \sum_{k=1}^{x}(a_{i,k} * m_{k,j}) \subseteq a_{i,j}$ as desired.

\par
This completes the proof.
$\square$ \\

}

\par Recall that if $R = \{0,1\}$ and $+$ and $*$ are given by logical OR and logical AND, respectively, then $\fdlcat{R}$ is called the category of finite sets with relations and is denoted $\REL$.  Therefore, the following corollary is a special case of Theorem 3.2.
\begin{cor}
The category of finite sets with relations, $\REL$, is dimension zero over any infinite field.
\end{cor}

\par
In fact, Bouc and Th\'{e}venaz \cite{correspond} have independently shown that $\REL$ is dimension zero over any field.  Additionally, they computed the irreducible representations of $\REL$ (Theorem 17.19 in \cite{correspond}) and showed that the fundamental correspondence functors (Definition 4.7 in \cite{correspond}) appear as subfunctors of a particular functor that arises from a lattice (Theorem 14.16 in \cite{correspond}).

\section{Acknowledgements}

\par
\indent \indent I would like to thank John Wiltshire-Gordon and David Speyer for working with me and for providing me with an invaluable undergraduate research experience.  I would also like to thank John Wiltshire-Gordon for providing mathematical insights and motivation which pushed me towards a proof of the main result in this paper.  I would like to thank the organizers of the University of Michigan REU program for giving me a wonderful opportunity to conduct mathematical research over the summer of 2015.  Finally, I would like to thank Stephen DeBacker for setting me up with the REU opportunity.  This research was partially supported by NSF Department of Undergraduate Education award 1347697 (REBUILD).

\pagebreak

\end{document}